\documentclass{amsart}

\begin{document}

\title[Evaluation of a new supply strategy for a fashion discounter]{Evaluation of a new supply strategy based on stochastic programming for a fashion discounter}
\author{Miriam Kie\ss ling, Tobias Kreisel, Sascha Kurz and J\"org Rambau}
\address{Business Mathematics, University of Bayreuth, 95440 Bayreuth, miriam.kiessling@uni-bayreuth.de, tobias.kreisel@uni-bayreuth.de, sascha.kurz@uni-bayreuth.de, j\"org.rambau@uni-bayreuth.de}

\maketitle

\begin{abstract}
  Fashion discounters face the problem of ordering the right
  amount of pieces in each size of a product.  The product is ordered
  in pre-packs containing a certain size-mix of a product. For this
  so-called lot-type design problem, a stochastic mixed integer linear
  programm was developed, in which price cuts serve as recourse action
  for oversupply. Our goal is to answer the question, whether the
  resulting supply strategy leads to a supply that is significantly
  more consistent with the demand for sizes compared to the original
  manual planning. Since the total profit is influenced by too many
  factors unrelated to sizes (like the popularity of the product, the
  weather or a changing economic situation), we suggest a comparison
  method which excludes many outer effects by construction. We apply
  the method to a real-world field study: The improvements in the size
  distributions of the supply are significant.
\end{abstract}

\section{Introduction}
\label{sec:1}
The lot-type design problem LDP seeks for an optimal set of lot-types
and a supply in terms of lots of such lot-types such that the
resulting supply of sizes matches the branch- and size-dependent
demand as closely as possible (see \cite{p_median} for details).  A
\emph{lot-type} is defined as a vector with a component for each size.
This component specifies the number of pieces of that size in the
lot. For example, for the sizes $S$, $M$, $L$, and $XL$, a lot of
lot-type $(1,2,3,1)$ contains one item of size~$S$, two of size~$M$,
three of size~$L$, and one of size~$XL$.

In reality, excess supply is compensated by cutting prices.
Therefore, we extended the model in \cite{p_median} by a model for
this recourse action. This resulted in the following stochastic mixed
integer linear programm, in the following denoted by SLDP:

 \begin{equation}
  \max  -\sum_{b \in B} \sum_{l \in L} \sum_{m \in M} c_{b,l,m} x_{b,l,m} - \sum_{i=1}^n \delta_i \cdot z_i + \mathbb{E}_\xi(Q(x,\xi))
\end{equation}
\begin{align}
  \text{s. t.}
  && 
  \sum_{l\in L}\sum_{m \in M} x_{b,l,m} &= 1 && \forall b\in B,\label{SOP1}\\
  &&
  \underline{I} \le \sum_{b\in B}\sum_{l\in L}\sum_{m \in M}
  m \cdot |l| \cdot x_{b,l,m} &\le \overline{I}, \label{SOP2}&&\\
  &&
  \sum_{l \in L} y_l & \le \sum_{i=1}^n z_i, \label{SOP3}&&\\
  &&
  \sum_{m \in M} x_{b,l,m} & \le y_l && \forall b\in B, 
  l\in L, \label{SOP4}
  \end{align}
  \begin{align} 
  &&
  z_i & \le z_{i-1} && \forall i \in \{2,\ldots,n\}, \label{SOP5}\\
  &&
  x_{b,l,m} & \in \{0,1\} && \forall b \in B, l \in L, m \in M, \label{SOP6}\\
  &&
  y_{l} & \in \{0,1\} && \forall l \in L, \label{SOP7}\\
  &&
  z_i & \in \{0,1\} && \forall i \in \{1,\ldots,n\}, \label{SOP8}
\end{align}
The meanings of the symbols are as follows: By~$L$ we denote the set
of possible lot-types which can be delivered to a branch from the
set~$B$ with multiplicity from the set~$M$. If Branch~$b$ gets
Lot-type~$l$ with Multiplicity~$m$ then the corresponding binary
variable $x_{b,l,m}$ takes value $1$, and $0$ otherwise. At most $n$
lot-types may be used. If we deliver Lot-type $l$ to Branch $b$ with
Multiplicity $m$ costs of $c_{b,l,m}$ arise. Every branch is supplied
by exactly one lot-type with one
multiplicity, see constraint~(\ref{SOP1}). The overall supply has to
be between the lower bound $\underline{I}$ and the upper bound
$\overline{I}$~(\ref{SOP2}). By $|l|$ we denote the overall number
of items in lot-type $l$. Using $i$ lot-types implies costs of
$\delta_i$. Constraint~(\ref{SOP5}) requires that also the costs for
using $i-1$, $i-2$, \dots, $1$ are added if $i$ different lot-types
are used. Due to Constraint~(\ref{SOP4}), $y_l$ indicates whether some
branch $b\in B$ is supplied by Lot-type $l$. Finally,~(\ref{SOP3})
links $y_l$ and~$z_i$.

To compensate for excess supply -- depending on a vector of random
variables $\xi$ describing the demand -- price cuts are possible. The
corresponding optimization problem seeking for optimal price-cut
strategies is denoted by $Q(x,\xi)$ where $x$ is the vector of all
variables $x_{b,l,m}$. \footnote{An approximation of $Q(x,\xi)$ can be
  computed by solving a mixed integer linear programm or a dynamic
  program.} The objective function is the revenue depending on
supply~$x$ and the random vector $\xi$. Our goal is to maximize the
total profit in expectation.

We performed a field study at our project partner using the SLDP. How
can we find out, whether the new method outperforms the traditional
manual planning? From an economic point of view, we should look at the
actual revenues. However, the data showed that actual revenues during
the field study are distorted too much by factors beyond our
control. For this purpose we developed two different comparison
methods which exclude such extern distortions and just reveil how well
the size-dependent demand is met. We examine significance by using the
Wilcoxon rank-sum test from statistics.

In the following we introduce our approaches and show results from the
field study.

\section{Comparison methods}
\label{sec:2}
We developed two different indicators for demand consistency of the
supply with sizes.

The first approach -- called Normalized Sales Rate Deviation
($\mathrm{NSRD}$) -- compares the sales rates per size against each
other in such a way that the popularity of the product itself has no
dominating influence on the result anymore.

The simple and course idea is to observe the sales at the
\emph{$50$\,\%-day} and look for each size at the fraction of supply
that has been sold so far.  The $50$\,\%-day is the first day in the
sales period where $50$~\% of the total supply of the product has been
sold.  We get an 
estimation of a normalized sales rate of each size relative to its
supply, independent of the popularity of the product. We call this
number the \emph{normalized sales rate estimate} of a size, denoted by
$\operatorname{NSR}(s)$.  

For example, if there is a product with supply $(10,20,20,10)$ 
, and the size-dependent sales numbers up to
the $50$\,\%-day are $(2,10,15,3)$, then the estimated sales rates are
$(0.2,0.5,0.75,0.3)$.  This indicates that the supply of M was spot-on
(i.e., relative sales in this size were the same as relative sales
in total), wheras the supplies of S and XL were too small and the supply
of L was too large in the considered branches.

\newcommand{\SDplus}{\ensuremath{\mathrm{SD}^{+}}} Based on these
observations, we are interested in how much the sales in a size
deviate from the overall $50$\,\%.  To this end, we estimate the
standard deviation of the observed normalized selling rates (relative
to the supply).  Because of multiple sales per day we can have an
average $\mathrm{NSR}$ different from $0.5$.  Therefore, the standard
deviation must be taken with respect to the sample average, which
is\footnote{We use the notation from~\cite{FPP}.}
\begin{equation}
\SDplus := \sqrt{\frac{1}{N - 1} \sum_{i=1}^{N}(r_{i} - \hat{r})^{2}},
\end{equation}
where $N$ is the sample size and $\hat{r}$ the sample average.  In our
case, we compute for $N$ sizes
\begin{equation}
  \label{eq:1}
  \mathrm{NSRD} := \sqrt{\frac{1}{N-1} \sum_{i=1}^{N}(\mathrm{NSR}(s_i) - \widehat{\mathrm{NSR}})^{2}},
\end{equation}
the \emph{normalized sales rate deviation} (for a single product).

The smaller the normalized sales rate deviation is, the more
consistent is the supply with the demand for sizes.

The second approach uses ideas from \cite{top_dog}. Let $P$ be the set
of products in a given commodity group, $S_p$ be the set of sizes for
a product $p\in P$, $b$ be a given branch, $s$ be a given size, and
$\theta_{b,s'}(p)$ be the first day at which Product $p$ is sold out
in Branch $b$ and Size $s'$ (here $\theta_{b,s'}(p)=\infty$ is
possible). With this we can define the Top-Dog-Count $W(b,s)$ as
\begin{equation}
   \left|\left\{p\in P\,\,\Big|\,\, 0=|\{s'\in S_p\mid
   \theta_{b,s'}(p)<\theta_{b,s}(p)\}|\right\}\right|
\end{equation}
and the Flop-Dog-Count $L(b,s)$ as
\begin{equation}
   \left|\left\{p\in P\,\,\Big|\,\, 0=|\{s'\in S_p\mid
   \theta_{b,s'}(p)>\theta_{b,s}(p)\}|\right\}\right|.
\end{equation}

We now consider the value $|W(b,s)-L(b,s)|$, in the following denoted
as \emph{top-dog-deviation of Size~$s$}, denoted by $\mathrm{TDD}(b,
s)$.  We assume that a value closer to zero corresponds to a better
supply strategy for $s$, because $\mathrm{TDD}(b, s)$ estimates the
difference of the probability that $s$ is sold out first and the
probability that $s$ is sold out last.  In order to exclude numerical
artefacts with too small integers, we restrict ourselves to samples
that consist of a (random) subset of
articles such that $W(b,s)+L(b,s)$ equals a given number. For each
such sample we compute $\mathrm{TDD}(b, s)$.

Since for each branch and each size we get an individual measurement,
the $\mathrm{TDD}$ is a finer measurement than the $\mathrm{NSRD}$.

Now, any improvement measured by the above indicators in the field
study could be a mere conincidence. To obtain statements about the
significance of the observations $\mathrm{NSRD}$ and $\mathrm{TDD}(b,
s)$, we apply the non-parametric \emph{Wilcoxon rank-sum test}
(see~\cite{wilcoxon}).
This method tests whether two sets of realizations stem from the
same distribution.  

Let us first sketch the Wilcoxon rank-sum test.  Given two sets of
realizations $A$ and $B$, the null hypothesis is that $A$ and $B$ stem
from the same distribution.  The alternative says that distribution of
$A$ is shifted to the left.

In our context the alternative means: there is less deviation among
normalized sales rates in the distribution of $A$ than there is in the
distribution of $B$.  Or, respectively, $A$ has smaller
$\mathrm{TDD}(b,s)$s than~$B$.  

The test is then conducted by sorting the values of both samples in
ascending order, thereby assigning a rank to each value. In the next
step the sum of the ranks for Family~$A$ is computed, which yields an
observed rank sum.

In order to estimate the probability that the distribution of~$A$ is
shifted to the left compared to the distribution of~$B$, we have to
consider the probability for getting rank sums smaller than or equal
to the observed rank sum for~$A$. If this probability lies below a
predetermined significance level (5\% or 10\%) we say that the
observed difference can not be explained by chance, and we can assume
that the $A$'s distribution is shifted to the left compared to the
distribution of~$B$.

\section{Comparison results}
\label{sec:3}
We employed the SLDP on results of a field study with ladieswear
shirts with sales periods from February to June 2011. For comparison
we took historical data from the same commodity group in a time period
from March until December 2006 -- at this time still all items were
supplied by manual planning.  We
consider 26 branches and four different sizes, namely S, M,
L and XL.

For the historical sample set only one lot-type was delivered, namely
$(1,2,2,1)$. In the field study our system provided the lot-types
$(1,1,1,1)$, $(1,1,2,2)$ and $(2,2,3,4)$. Now we would like to find
out whether SLDP led to supplies that were significantly more
consistent with the demand than the supplies suggested by manual
planning.

At first we have a look at the normalized sales rate
deviation. Aggregating over the tested branches and rounding to two
decimal figures leads us to the observations for $\mathrm{NSRD}$ with
the corresponding ranks in braces in Table \ref{table:sr}.

\begin{table}[htp]
\begin{center}
\tiny
\begin{tabular}{cccccccc}
 \textbf{SLDP} 
& 8.58 (4)  & 16.92 (18) & 13.49 (13) & 17.32 (21) & 17.29 (20) & 14.37 (14) & 12.83 (12) \\
& 9.85 (5) & 6.61 (3) & 16.40 (17) & & & \\
&&&&&&& \\
 \textbf{manual planning} 
& 18.15 (24) & 21.89 (28) & 22.71 (31) & 12.44 (9) & 19.93 (26) & 22.98 (33) & 27.46 (35) \\
& 21.51 (27) & 22.79 (32) & 15.78 (15) & 12.49 (10) & 16.09 (16) & 22.62 (30) & 19.09 (25) \\ 
& 17.06 (19) & 10.52 (6) & 17.53 (22) & 11.66 (7) & 12.35 (8) & 5.09 (2) & 17.65 (23) \\
& 3.37 (1) & 22.08 (29) & 12.75 (11) & 23.92 (34) & 27.52 (36) & 42.32 (38) & 36.80 (37) 
\end{tabular}
\caption{standard deviation of selling rates ($\mathrm{NSRD}$) for
  tested articles of ladieswear with ranks in braces} 
\label{table:sr}
\end{center}
\end{table}

We get a mean of 13.36 for the SLDP and 19.16 for the manual
planning. So in average we get smaller deviations among the normalized 
sales rates per size for SLDP.  Thus, we conclude that the systematic
error in size distribution is smaller. To test significance, we apply
the Wilcoxon rank-sum test with a predefined significance level of
5\%. The null hypothesis is that the distributions of $\mathrm{NSRD}$
for SLDP and the manual planning are identical, the alternative that
for SLDP the distribution is shifted to the left. The test yields a
rank
sum of 127 for the SLDP. The probability of a rank sum lower than
or equal to 127 is approximately 1.17\%, which is below the significance
level. In other words: The probability that the observed improvements
are the result of coincidence is well below the significance level 5\%, and the
observations are significant.

As a next step we check the measurements
$\mathrm{TDD}(b,s)$. Analogous to the selling rates, we compute mean
values over all products, branches, and sizes, denoted by
$\widehat{\mathrm{TDD}}$. For the SLDP, $\widehat{\mathrm{TDD}}$ is 2,
whereas it is 3 for the manual planning. Again, we perform the
Wilcoxon rank-sum test with a significance level of 5\%. We get the
following result: With 208 observations (104 for each sample set) and a
rank sum of 8766 for the SLDP method, the probability for a shift to
the left of the distribution of $\widehat{\mathrm{TDD}}$ for SLDP
method is less than $10^{-6}$\%.  This also is well below the
significance level: The observed improvements are significant.

\section{Conclusion}
\label{sec:4}

We presented methods to compare the impact of different supply
strategies for a fashion discounter in a real-world field study.  The
normalized sales rate deviation among sizes measures how evenly a
product sells in the various sizes.  The top dog deviation measures to
what extent a size is sold out first more often than last, or vice
versa.  Both measures showed that our new supply strategy based on the
stochastic lot design problem can significantly improve the demand
consistency of the supply with respect to sizes, where significance is
checked by the Wilcoxon rank sum test.  Since the Wilcoxon text is
robust and uses no assumptions on the distributions, we are confident
that this result is practically relevant.


\begin{thebibliography}{99.}
\bibitem{p_median}
C.~Gaul, S.~Kurz and J.~Rambau: On the lot-type design problem, Optimization Methods and Software Volume 25, Issue 2 (2010), p. 217-227.

\bibitem{wilcoxon}
H.~B\"uning and G.~Trenkler: Nichtparametrische statistische Methoden, Gruyter, 1994

\bibitem{FPP}
R.~Purves, D.~Freedman and R.~Pisani: Statistics, WW Norton \& Co,
  1998.

\bibitem{paper_rio}
S.~Kurz and J.~Rambau: Demand forecasting for companies with many
  branches, low sales numbers per product, and non-recurring orderings,
  Proceedings of the Seventh International Conference on Intelligent Systems
  Design and Applications, 22-24.10.2007, Rio de Janeiro, Brazil, 2007,
  pp.~196--201.

\bibitem{top_dog}
S.~Kurz and J. Rambau, J.~Schl\"uchtermann, and R.~Wolf.: The {T}op-{D}og
  {I}ndex: A {N}ew {M}easurement for the {D}emand {C}onsistency of the {S}ize
  {D}istribution in {P}re-{P}ack {O}rders for a {F}ashion {D}iscounter with
  {M}any {S}mall {B}ranches,  (submitted).
\end{thebibliography}
\end{document}